
\input amstex
\input amsppt.sty
\hsize 30pc
\vsize 47pc
\def\cit#1#2{\ifx#1!\cite{#2}\else#2\fi} 
\def\idx{}               
\def\ign#1{}             

\redefine\o{\circ}
\define\X{\frak X}
\define\al{\alpha}
\define\be{\beta}
\define\ga{\gamma}

\define\ep{\varepsilon}

\define\la{\lambda}

\define\si{\sigma}

\define\ph{\varphi}

\define\om{\omega}

\define\Om{\Omega}
\redefine\i{^{-1}}
\define\row#1#2#3{#1_{#2},\ldots,#1_{#3}}

\define\x{\times}
\define\Vect{\operatorname{Vect}}
\define\Diff{\operatorname{Diff}}
\define\Fl{\operatorname{Fl}}
\redefine\div{\operatorname{div}}
\redefine\L{\Cal L}
\define\grad{\operatorname{grad}}
\def\today{\ifcase\month\or
 January\or February\or March\or April\or May\or June\or
 July\or August\or September\or October\or November\or December\fi
 \space\number\day, \number\year}
\topmatter
\title  $n$-transitivity of certain diffeomorphism groups
\endtitle
\author  Peter W. Michor \\
Cornelia Vizman \endauthor
\affil
Erwin Schr\"odinger International Institute of Mathematical Physics,
Wien, Austria
\endaffil
\address
P. Michor: Institut f\"ur Mathematik, Universit\"at Wien,
Strudlhofgasse 4, A-1090 Wien, Austria
\endaddress
\email Peter.Michor\@esi.ac.at \endemail
\address
C. Vizman: University of Timisoara, Mathematics Departement,
Bul.V.Parvan 4, 1900-Timisoara, Romania
\endaddress
\dedicatory \enddedicatory
\date {\today} \enddate
\thanks \endthanks
\abstract It is shown that some groups of diffeomorphisms of a
manifold act $n$-transitively for each finite $n$.
\endabstract
\endtopmatter

\document


Let $M$ be a connected smooth manifold of
dimension $\dim M\ge2$.
We say that a subgroup $G$ of the group $\Diff(M)$ of all smooth
diffeomorphisms acts \idx{\it $n$-transitively on $M$},
if for any two ordered sets of $n$ different points
$(x_1,\dots,x_n)$ and $(y_1,\dots,y_n)$ in $M$ there is a smooth
diffeomorphism $f\in G$ such that $f(x_i)=y_i$ for
each $i$.

\proclaim{Theorem} Let $M$ be a connected smooth (or real
analytic) manifold of dimension $\dim M\ge2$.
Then the following subgroups of the group $\Diff(M)$ of all smooth
diffeomorphisms with compact support act $n$-transitively on $M$,
for each finite $n$:
\roster
\item The group $\Diff_c(M)$ of all smooth diffeomorphisms with
       compact support.
\item The group $\Diff^\om(M)$ of all real analytic diffeomorphisms.
\item If $(M,\si)$ is a symplectic manifold, the group $\Diff_c(M,\si)$
       of all symplectic diffeomorphisms with compact support, and
       even the subgroup of all globally Hamiltonian
       symplectomorphisms.
\item If $(M,\si)$ is a real analytic symplectic manifold, the group
       $\Diff^\om(M,\si)$
       of all real analytic symplectic diffeomorphisms, and even the
       subgroup of all globally Hamiltonian real analytic
       symplectomorphisms.
\item If $(M,\mu)$ is a manifold with a smooth volume density, the
       group $\Diff_c(M,\mu)$ of all volume preserving
       diffeomorphisms with compact support.
\item If $(M,\mu)$ is a manifold with a real analytic volume density, the
       group $\Diff^\om(M,\mu)$ of all real analytic volume preserving
       diffeomorphisms.
\item If $(M,\al)$ is a contact manifold, the group $\Diff_c(M,\al)$
       of all contact diffeomorphisms with compact support.
\item If $(M,\al)$ is a real analytic contact manifold, the group
       $\Diff^\om(M,\al)$ of all real analytic contact diffeomorphisms.
\endroster
\endproclaim

Result \therosteritem1 is folklore, the first trace is in \cit!{8}.
The results \therosteritem3, \therosteritem5, and \therosteritem7 are
due to \cit!{3} for 1-transitivity, and to \cit!{1} in the general
case. Result \therosteritem2 is from \cit!{7}. We shall give here a
short uniform proof, following an argument from \cit!{7}. That this
argument suffices to prove all results was noted by the referee, many
thanks to him.

\demo{Proof}
Let us fix a finite $n\in \Bbb N$.
Let $M^{(n)}$ denote the open submanifold of all $n$-tuples
$(x_1,\dots,x_n)\in M^{n}$ of pairwise  distinct points.
Since $M$ is connected and of dimension $\ge 2$, each $M^{(n)}$ is
connected.

The group $\Diff(M)$ acts on $M^{(n)}$ by the diagonal action, and we
have to show, that any of the subgroups $G$ described above acts
transitively. We shall show below that for each $G$ the $G$-orbit
through any $n$-tuple $(x_1,\dots,x_n)\in M^{(n)}$ contains an open
neighborhood of $(x_1,\dots,x_n)$ in $M^{(n)}$, thus any orbit is
open; since $M^{(n)}$ is connected, there can then be only one orbit.
\qed\enddemo

\proclaim{Lemma}
Let $M$ be a real analytic manifold.
Then for any real analytic vector bundle $E\to M$ the space $C^\om(E)$ of
real analytic sections of $E$ is dense in the space $C^\infty(E)$ of
smooth sections.
In particular
the space $\X^\om(M)$ of real analytic vector fields is dense in the
space $\X(M)$ of smooth vector fields, in the Whitney
$C^\infty$-topology.
\endproclaim

\demo{Proof}
For functions instead of sections this is \cit!{2}, proposition 8.
Using results from \cit!{2} it can easily be extended to sections,
as is done in \cit!{6}, 7.5.
\qed\enddemo

\demo{\bf The cases \therosteritem2 and \therosteritem1}
We choose a complete Riemannian metric $g$ on $M$ and we let
$(Y_{ij})_{j=1}^m$ be an orthonormal basis of $T_{x_i}M$ with respect
to $g$, for all $i$. Then we choose real analytic vector fields
$X_k$ for $1\le k\le N=nm$ which satisfy the
following conditions:
$$\alignat2
&|X_k(x_i)-Y_{ij}|_g<\ep&\quad&\text{ for }k=(i-1)m+j,\\
&|X_k(x_i)|_g<\ep&\quad&\text{ for all }k\notin [(i-1)m+1, im],\tag9\\
&|X_k(x)|_g<2 &&\text{ for all }x\in M \text{ and all }k.
\endalignat$$
Since these conditions describe a Whitney $C^0$ open set, such vector
fields exist by the lemma.
The fields are bounded with respect to a complete Riemannian
metric, so they have complete real analytic flows $\Fl^{X_k}$,
see e\.g\. \cit!{4}.
We consider the real analytic mapping
$$\gather
f:\Bbb R^N\to M^{(n)}\\
f(t_1,\dots,t_N):=
\pmatrix (\Fl^{X_1}_{t_1}\o\dots\o\Fl^{X_N}_{t_N})(x_1)\\
	\dotsc \\ (\Fl^{X_1}_{t_1}\o\dots\o\Fl^{X_N}_{t_N})(x_n)
\endpmatrix
\endgather$$
which has values in the $\Diff^\om(M)$-orbit through $(x_1,\dots,x_n)$.
To get the tangent mapping at 0 of $f$ we consider the partial
derivatives
$$
\frac{\partial}{\partial t_k}|_0 f(0,\dots,0,t_k,0,\dots,0) =
     (X_k(x_1),\dots,X_k(x_n)).
$$
If $\ep>0$ is small enough, this is near an orthonormal basis of
$T_{(x_1,\dots,x_n)}M^{(n)}$ with respect to the product metric
$g\x \dots \x g$. So $T_0f$ is invertible and the image of $f$
contains thus an open subset.

In case \therosteritem1, we can choose smooth vector
fields $X_k$ with compact support which satisfy conditions \thetag9.
\qed\enddemo

For the remaining cases we just indicate the changes which are
necessary in this proof.

\demo{\bf The cases \therosteritem4 and \therosteritem3}
Let $(M,\si)$ be a connected  real analytic
symplectic smooth manifold of dimension $m\ge 2$.
We choose real analytic functions
$f_k$ for $1\le k\le N=nm$ whose Hamiltonian vector fields
$X_k=\grad^\si(f_k)$ satisfy conditions \thetag9.
Since these conditions describe Whitney $C^1$ open subsets, such
functions exist by \cit!{2}, proposition 8. Now we may finish the
proof as above.
\qed\enddemo

\subhead Contact manifolds \endsubhead
Let $M$ be a smooth manifold of dimension $m=2n+1\ge 3$. A
\idx{\it contact form } on $M$ is a 1-form $\al\in\Om^1(M)$ such that
$\al\wedge (d\al)^n\in\Om^{2n+1}(M)$ is nowhere zero. This is
sometimes called an \idx{\it exact} contact structure. The pair
$(M,\al)$ is called a \idx{\it contact manifold} (see
\cit!{5}).
The \idx{\it contact vector field} $X_\al\in\X(M)$ is the unique
vector field satisfying $i_{X_\al}\al=1$ and $i_{X_\al}d\al=0$.

A diffeomorphism $f\in\Diff(M)$ with $f^*\al=\la_f.\al$ for a nowhere
vanishing function $\la_f\in C^\infty(M,\Bbb R\setminus 0)$ is called
a \idx{\it contact diffeomorphism}. Note that then
$\la_f = i_{X_\al}(\la_f.\al) = i_{X_\al}f^*\al =
f^*(i_{(f\i)^*X_\al}\al) = f^*(i_{f_*X_\al}\al)$.
The group of all contact diffeomorphisms will be
denoted by $\Diff(M,\al)$.

A vector field $X\in\X(M)$ is called a contact vector field if
$\L_X\al=\mu_X.\al$ for a smooth function
$\mu_X\in C^\infty(M,\Bbb R)$.
The linear space of all contact vector fields will be denoted by
$\X_\al(M)$ and it is clearly a Lie algebra. Contraction with $\al$ is a
linear mapping again denoted by
$\al:\X_\al(M) \to C^\infty(M,\Bbb R)$.
It is bijective since we may apply $i_{X_\al}$ to the equation
$\L_X\al = i_X\,d\al + d \al(X) = \mu_X.\al$
and get $0+i_{X_\al}d\al(X) = \mu_X$; but the equation uniquely
determines $X$ from $\al(X)$ and $\mu_X$. The inverse
$f\mapsto \grad^\al(f)$ of $\al:\X_\al(M)\to C^\infty(M,\Bbb R)$ is a linear
differential operator of order 1.

\demo{\bf The cases \therosteritem8 and \therosteritem7}
Let $(M,\al)$ be a connected real analytic contact
manifold of dimension $m\ge2$.
We choose real analytic functions
$f_k$ for $1\le k\le N=nm$ such that their contact vector fields
$X_k=\grad^\al(f_k)$ satisfy conditions \thetag9.
Since these conditions describe Whitney $C^1$ open subsets, such
functions exist by \cit!{2}, proposition 8. Now we may finish the
proof as above.
\qed\enddemo

\demo{\bf The cases \therosteritem6 and \therosteritem5}
Let $(M,\mu)$ be a connected real analytic
manifold of dimension $m\ge 2$ with a real analytic positive volume density.
We can find a real analytic Riemannian metric $\ga $ on $M$ whose volume form
is $\mu$.
Then the divergence of a vector field $X\in\Vect(M)$
is $\div X=*d*X^\flat$, where $X^\flat=\ga (X)\in\Om^1(M)$  (here we
view $\ga :TM\to T^*M$)
and $*$ is the Hodge star operator of $\ga$.
We also choose a complete Riemannian metric $g$.

First we assume that $M$ is orientable.
We choose real analytic $(m-2)$-forms $\be_k$ for $1\le k\le N=nm$
such that the vector fields $X_k=(-1)^{m+1}\ga \i*d\be_k$ satisfy
conditions \thetag9.
Since these conditions describe Whitney $C^1$ open subsets, such
$(m-2)$-forms exist by the lemma. The real analytic vector fields
$X_k$ are then divergence free since
$\div X_k=*d*\ga X_k=*dd\be_k=0$.
Now we may finish the proof as usual.

For non-orientable $M$, we let $\pi:\tilde M\to M$
be the real analytic connected oriented double cover of $M$, and let
$\ph:\tilde M\to \tilde M$ be the real analytic involutive covering
map.
We let $\pi\i(x_i)=\{x^1_i,x^2_i\}$, and we
pull back both metrics to $\tilde M$, so $\tilde \ga :=\pi^*\ga $ and
$\tilde g:=\pi^*g$.
We choose real analytic $(m-2)$-forms
$\be_k\in\Om^{m-2}(\tilde M)$ for $1\le k\le N=nm$ whose vector fields
$X_{\be_k}=(-1)^{m+1}\tilde \ga \i*d\be_k$ satisfy the
following conditions, where we put
$Y^p_{ij}:=T_{x^p_{ij}}\pi\i.Y_{ij}$ for $p=1,2$:
$$\alignat2
&|X_{\be_k}(x^p_i)- Y^p_{ij}|_{\tilde g}<\ep
     &\quad&\text{ for }k=(i-1)m+j,p=1,2,\\
&|X_{\be_k}(x^p_i)|_{\tilde g}<\ep
     &\quad&\text{ for all }k\notin [(i-1)m+1, im],p=1,2,  \tag{10}\\
&|X_{\be_k}|_{\tilde g}<2
     &&\text{ for all }x\in \tilde M \text{ and all }k.
\endalignat$$
Since these conditions describe Whitney $C^1$ open subsets, such
$(m-2)$-forms exist by the lemma. Then the vector fields
$\frac12(X_{\be_k} + \ph_*X_{\be_k})$ still satisfy the conditions
\thetag{10}, are still divergence free and induce divergence free
vector fields $Z_{\be_k}\in\X(M)$ which satisfy the conditions
\thetag9 on $M$ as in the oriented case, and we may finish the proof
as above.
\qed\enddemo

\enddocument